\documentclass[12pt,twoside]{article}
\usepackage{amsmath,amsbsy,amsfonts}
\def\dis{\displaystyle}

\newtheorem{lemma}{Lemma}[section]
\newtheorem{proposition}{Proposition}[section]
\newtheorem{theorem}{Theorem}[section]
\newtheorem{corollary}{Corollary}[section]

\newcommand{\fim}{\hfill\rule{2mm}{2mm}}

\def\dis{\displaystyle}

\let\Section=\section
\def\section{\setcounter{equation}{0}\Section}

\begin{document}

\title{\Large\sf{Multiplicity of solutions for a NLS equations with magnetic fields in $\mathbb{R}^{N}$}}

\date{}

\author{\sf{Claudianor O. Alves} \thanks{Supported by INCT-MAT, PROCAD,NPq/Brazil
620150/2008-4 and 303080/2009-4}\\
{\small \textit{Unidade Acad\^emica de Matem\'atica e Estat\'{\i}stica}}\\
{\small \textit{Universidade Federal de Campina Grande}}\\
{\small \textit{58429-900, Campina Grande - PB - Brazil}}\\
{\small \textit{e-mail address: coalves@dme.ufcg.edu.br}} \\
\\
\vspace{1mm} \sf{Giovany M. Figueiredo} \thanks{Research partially
supported by CNPq/PQ
301242/2011-9 and 200237/2012-8.}\\
{\small \textit{Faculdade de Matem\'atica}}\\
{\small \textit{Universidade Federal do Par\'a}}\\
{\small \textit{66075-110, Bel\'em - PA - Brazil}}\\
{\small \textit{e-mail address: giovany@ufpa.br}}} \maketitle

\begin{abstract}
\noindent We investigate the multiplicity of nontrivial weak
solutions for a class of complex equations.  This class of problems
are related with the existence of solitary waves for a nonlinear
Sch\"{o}dinger equation. The main result is established by using
minimax methods and Lusternik-Schnirelman theory of critical points.
\end{abstract} \maketitle

\bigskip



\vskip.2cm

\section{Introduction}

In this paper, we establish existence and multiplicity of nontrivial
weak solutions for the following class of nonlinear Schr\"{o}dinger
equations:
$$
\left\{
 \begin{array}{ll}
\left( \displaystyle\frac{\epsilon}{i}\nabla - A(z)\right)^2 u +
(\lambda W(z)+1)u = f(|u|^2)u,&~~z \in \mathbb{R}^N, \ \
\\
u \in H^1(\mathbb{R}^N,\mathbb{C}),
\end{array}
       \right.
           \eqno{(P_{\epsilon,\lambda})}
$$
where  $\varepsilon , \lambda $ are positive parameters and $W,f$ are continuous functions satisfying some technical conditions.

This class of problem is related with the existence of solitary
waves, namely solutions of the form $\psi(x,t) := e^{-i
\frac{E}{\epsilon}t}u(x)$, with $E \in \mathbb{R}$, for a nonlinear
Schr\"odinger equation like
$$
i \epsilon \dfrac{ \partial\psi}{ \partial t } = \left( \dfrac{\epsilon}{i}
\nabla - A(z) \right)^2 \psi +U(z)\psi - f(|\psi|^2)\psi,~~z \in
\mathbb{R}^{N}, \eqno{(NLS)}
$$
where $t >0$, $N \geq 2$, $\epsilon$ is the Planck constant and $A$ is a magnetic potential associated to a given
magnetic $B$, $U(z)$ is a real electric potential and the nonlinear
term $f$ is a superlinear function. A direct computation shows
that $\psi$ is a solitary wave for $(NLS)$ if, and only if, $u$ is
a solution of the following problem
\begin{equation} \label{pe}
\left( \dfrac{\epsilon}{i}\nabla - A(z)\right)^2 u + (\lambda W(z)+1)u =
f(|u|^2)u,~~\mbox{ in } \mathbb{R}^{N},
\end{equation}
where $\lambda W(z)+1 = U(z)-E$. It is important  to investigate the
existence and the shape of such solutions in the semiclassical
limit, namely, as $\epsilon \to 0^+$. The importance of this study relies
on the fact that the transition from Quantum Mechanics to
Classical Mechanics can be formally performed by sending the
Planck constant to zero.

At the last years, a lot of papers have considered some classes of problem related to problem $(P_{\epsilon,\lambda})$, for the case where the magnetic field is nontrivial, that is, $A \not= 0$, motivated by  a seminal paper due to Esteban and Lions \cite{EL}, see for example, Alves, Figueiredo \& Furtado \cite{AlvFigFur1,AlvFigFur2}, Cingolani \& Secchi \cite{CinSec1,CinSec2}, Cingolani, Jeanjean \& Secchi \cite{CinJeaSec}, Chabrowski \& Szulkin \cite{ChaSzu}, Ding \& Liu \cite{DX}, Ding \& Wang \cite{DZ}, Kurata \cite{Kurata}, Liang \& Zhang \cite{LZ}, Tang \cite{Tan} and their references.

In \cite{BW1}, Barstch \& Wang treated the problem (\ref{pe}) considering the case where $A=0,f(t)=|t|^{\frac{q-2}{2}}$ and $W$ satisfying :
\begin{itemize}
\item[$(W_{1})$]\,\, $W \in C(\mathbb{R}^{N}, \mathbb{R})$, $W(x) \geq
0$ for all $x \in \mathbb{R}^{N}$ and $\Omega={\rm{int}}\,
W^{-1}(0)$ is a nonempty bounded open set with smooth boundary
$\partial \Omega$ and $\overline{\Omega} = W^{-1}(0)$.

\item[$(W_{2})$] There exists $K_{0}>0$ such that
$$
\mu \Big(\Big\{ x \in \mathbb{R}^{N}: \, W(x) \leq K_{0}
\Big\}\Big)< \infty.
$$
\end{itemize}
In \cite{BW1}, it is proved the existence of a least energy solution
$u_\lambda$ of $(P)$ and for any sequence $\lambda_n \to \infty$ has
a subsequence such that   $u_{\lambda_n}$ converges strongly in
$H^1(\mathbb{R}^N)$ along the subsequence to a least energy solution
of the following limit problem
$$
\left\{
\begin{array}{l}
-\Delta{u}+ u=u^q, \,\,\, \mbox{in}\,\, \Omega\\
\mbox{}\\
u(x)>0 \,\,\, \mbox{in} \,\, \Omega\,\,\ \mbox{and}\,\,\, u=0 \,\,\,
\mbox{on} \,\, \partial\Omega
\end{array}
\right.
$$
Moreover, it is stated that there exist $1< q_0 < \frac{N+2}{N-2}$
and a function  $\Lambda:(q_0, \frac{N+2}{N-2}) \rightarrow
\mathbb{R}$ such that $(P)$ has at least ${\rm{cat}}(\Omega)$
solutions for any $\lambda \geq \Lambda(q)$.

Motivated by results showed in \cite{BW1}, Alves \& Soares in \cite{AlvesSoares} have considered the existence and multiplicity of solutions for the following class of quasilinear problem
$$
\left\{
\begin{array}{l}
-\epsilon^{p}\Delta_{p}{u}+(\lambda W(x)+1)|u|^{p-2}u=f(u), \,\,\, \mbox{in} \,\,\mathbb{R}^{N}\\
\mbox{}\\
u(x)>0 \,\,\, \mbox{in} \,\, \mathbb{R}^{N}
\end{array}
\right.
\leqno{(P_{\lambda, \epsilon})}
$$
where $\epsilon$ and $\lambda$ are positive parameters, $\Delta_{p}$
is the p-Laplacian operator and $2 \leq p <N$.  The assumptions on
$W$ are essentially those assumed in \cite{BW1}, namely  $(W_2)$ and
the following version of $(W_1)$:

\begin{itemize}

\item[$(W_{1}')$] \,\, $W \in C^{1}(\mathbb{R}^{N}, \mathbb{R}), W(x)
\geq 0$ for all $x \in \mathbb{R}^{N}$ and $\Omega={\rm{int}}
W^{-1}(0)$ is a nonempty bounded open set with smooth boundary
$\partial \Omega$ and $0 \in \Omega$. Moreover,
$W^{-1}(0)=\overline{\Omega} \cup D$ where $D$ is a set of measure
zero.
\end{itemize}
An important point related to conditions $(W_{1}')-(W_{2})$ is the fact that potential $V(x)=\lambda W(x)+1$ does not verifies the condition
$$
\displaystyle\liminf_{|x|\rightarrow \infty}V(x) >
\displaystyle\inf_{x \in \mathbb{R}^{N}} V(x)=V_0>0, \eqno{(R)}
$$
which has been introduced by Rabinowitz \cite{Rabinowitz}. By assuming some technical conditions on $f$, the main result proved in \cite{AlvesSoares} claims that if $(W_{1}')-(W_{2})$ hold, then there exists $\epsilon^{*}>0$ such that for any $\epsilon \in (0,
\epsilon^{*})$ there exists $\lambda^{*}(\epsilon)>0$ such that
$(P_{\lambda,\epsilon})$ has at least $cat(\Omega)$ solutions for
any $\lambda \geq \lambda^{*}(\epsilon)$.

The main proposal of the present is to establish the same type of result found in \cite{AlvesSoares} for problem $(P_{\epsilon,\lambda})$. To this end, we assume that the nonlinearity $f:\mathbb{R}\rightarrow\mathbb{R}$
is of class $C^1$ and satisfies the following conditions:

\begin{enumerate}
\item[$(f_1)$]
$f(s)=0 $ for all $s\leq 0$ and  $f(s)=o(|s|)$ at the origin;
\item[$(f_2)$] $\displaystyle \lim_{|s|\rightarrow
\infty}f(s)|s|^{(-q+2)/2}=0$ for some $q\in (2,2^{*})$, where
$2^*=2N/(N-2)$;
\item[$(f_3)$] There exists $\theta > 2$ such that
$0<\frac{\theta}{2} F(s)\leq sf(s)$ for all $s>0$.
\item[$(f_4)$] The function $s\rightarrow f(s)$ is
increasing for $s>0$.
\end{enumerate}

A typical example of a function satisfying the conditions
$(f_1)-(f_4)$ is given by $f(s)=s^{q_{1}}$ for $s\geq 0$, with
$1<q_{1}<q-1$, and $f(s)=0$ for $s<0$.\\

Before to state our main result, we recall that if $Y$ is a closed set of a topological space $X$, we denote the Lusternik-Schnirelmann
category of $Y$ in $X$ by $cat_X (Y)$, which is the least number of closed and contractible sets in $X$ that cover $Y$. Hereafter, $cat \, X$
denotes $cat_X (X)$.

\vspace{0.5 cm}

Our main result is the following:

\begin{theorem} \label{T1}
Suppose that $A \in C(\mathbb{R}^{N},\mathbb{R}^{N})$ is bounded,
$(W_{1}')-(W_{2})$ and $(f_{1})-(f_{4})$ hold. Then there exists
$\epsilon^{*}>0$ such that for any $\epsilon \in (0, \epsilon^{*})$
there exists $\lambda^{*}(\epsilon)>0$ such that
$(P_{\epsilon,\lambda})$ has at least $cat(\Omega)$ solutions for
any $\lambda \geq \lambda^{*}(\epsilon)$.
\end{theorem}

In the proof of Theorem \ref{T1},  we will use variational methods and prove important estimates involving some minimax levels. Here, since we are with the presence of the a magnetic field, some estimates involving the case $A=0$ cannot be repeat and new estimates are necessary to get the results, see for example the Sections 4 and 5, where key estimates involving the some minimax levels and barycenter were made.

The plan of this paper is as follows. In Section 2, we recall some properties involving the function spaces and the energy functional associated with problem $(P_{\epsilon,\lambda})$. In Section 3, we study the behavior of the Palais-Smale sequence involving the energy functional. In Section 4,  we show some estimates involving some minimax levels,  and finally in  Section 5, we prove Theorem \ref{T1}.

\section{Variational Framework}

By the change of variables $z\mapsto \epsilon x$, we can see that
$(P_{\varepsilon, \lambda})$ is equivalent to
$$
 \left\{
 \begin{array}{ll}
\left( \dfrac{1}{i}\nabla - A(\epsilon x)\right)^2 u +(\lambda
W(\epsilon x)+1)u = f(|u|^2)u, \ \ x \in \mathbb{R}^N,
\\
u \in H^1(\mathbb{R}^N,\mathbb{C}).
\end{array}
       \right.
           \eqno{(D_{\epsilon,\lambda})}
$$

For each $\epsilon,\lambda>0$, we set the Hilbert space
$$
H_{\epsilon,\lambda}=\biggl\{u \in H^{1}(\mathbb{R}^{N},\mathbb{C}):
\displaystyle\int_{\mathbb{R}^{N}}W(\epsilon x)|u|^{2}<
\infty\bigg\}
$$
under the scalar product
$$
\langle u,v \rangle_{\epsilon,\lambda} := \mbox{Re} \left(
\displaystyle\int_{\mathbb{R}^{N}} \nabla_{\epsilon} u
\overline{\nabla_{\epsilon} v}+(\lambda W(\epsilon
x)+1)u\overline{v} \right),
$$
where $\mbox{Re}(w)$ denotes the real part of $w \in \mathbb{C}$,
$\overline{w}$ is its conjugated, \linebreak $\nabla_{\epsilon} u :=
(D_{1}^{\epsilon}u,D_{2}^{\epsilon}u, ..., D_{N}^{\epsilon}u )$ and
$D_{j}^{\epsilon}:=i^{-1}\partial_{j}-A_j(\epsilon x )$, for
$j=1,\ldots,N$. The norm induced by this inner product is given by
$$
\|u\|_{\epsilon,\lambda}^{2}=
\left(\displaystyle\int_{\mathbb{R}^{N}} |\nabla_{\epsilon}
u|^2+(\lambda W(\epsilon x)+1)|u|^2\right)^{1/2}.
$$

In this work, the usual Hilbert space
$H^{1}(\mathbb{R}^{N},\mathbb{C})$ is provided with the inner
product
$$
\langle u,v \rangle_{\epsilon} := \mbox{Re} \left(
\displaystyle\int_{\mathbb{R}^{N}} \nabla_{\epsilon} u
\overline{\nabla_{\epsilon} v}+u\overline{v} \right).
$$

The norm induced by this inner product is given by
$$
\|u\|_{\epsilon}= \left(\displaystyle\int_{\mathbb{R}^{N}}
|\nabla_{\epsilon} u|^2+|u|^2\right)^{1/2}.
$$

As proved by Esteban and Lions in \cite[Section II]{EL}, for any $u
\in H_{\epsilon,\lambda}$ there holds
\begin{equation} \label{diamagnetica}
 |\nabla |u|(x)|= \left|\mbox{Re}\left(\nabla u
\dis\frac{\overline{u}}{|u|}\right)\right|=\left|\mbox{Re}\biggl((\nabla
u - iA_{\epsilon}u)\dis\frac{\overline{u}}{|u|}\biggl)\right|\leq
|\nabla_{\epsilon} u(x)|.
\end{equation}
The above expression is the so called diamagnetic inequality. It
follows from it that, if $u \in H_{\epsilon,\lambda}$, then $|u| \in
H^1(\mathbb{R}^N,\mathbb{R})$. Moreover, the embedding \linebreak
$H_{\epsilon,\lambda} \hookrightarrow L^q(\mathbb{R}^N,\mathbb{R})$
is continuous for each $2 \leq q \leq 2^*$ and, for each bounded set
$\Lambda \subset \mathbb{R}^N$ and $2\leq q<2^*$, the embedding
below is compact
\begin{equation}\label{imersaocompacta}
H_{\epsilon,\lambda} \hookrightarrow L^q(\Lambda,\mathbb{R}).
\end{equation}

 We say that a
function $u\in H_{\epsilon,\lambda}$ is a weak solution of the
problem $(P_{\epsilon})$ if
$$
\mbox{Re}\biggl( \dis\int_{\mathbb{R}^{N}} \nabla_{\epsilon} u
\overline{\nabla v_{\epsilon}}+(\lambda W(\epsilon
x)+1)u\overline{v}-f(|u|^{2})u\overline{v}\biggl)=0,~~\mbox{for each
}v \in H_{\epsilon,\lambda}.
$$
In view of $(f_2)$ and $(f_3)$, we have that the associated
functional $I_{\epsilon,\lambda} : H_{\epsilon,\lambda} \to
\mathbb{R}$ given by
$$
I_{\epsilon,\lambda}(u) := \dis\frac{1}{2}\int |\nabla_{\epsilon}
u|^2 +
          \frac{1}{2}\dis\int_{\mathbb{R}^{N}}( \lambda W(\epsilon x)+1)|u|^2 -\dis\frac{1}{2}
          \dis\int F(|u|^{2})
$$
is well defined. Moreover, $I_{\epsilon,\lambda}  \in
C^1(H_{\epsilon,\lambda}),\mathbb{R})$ with the following derivative
$$
I'_{\epsilon,\lambda}(u)v=\mbox{Re}\biggl( \dis\int_{\mathbb{R}^{N}}
\nabla_{\epsilon} u \overline{\nabla_{\epsilon} v}+(\lambda
W(\epsilon x)+1)u\overline{v}-f(|u|^{2})u\overline{v}\biggl) \,\,\, \forall u,v \in H_{\epsilon,\lambda}.
$$
Hence, the weak solutions of $(D_{\epsilon,\lambda})$ are  precisely
the critical points of $I_{\epsilon,\lambda}$.

\section{The Palais-Smale condition}

In this Section, the main goal is to show that functional $I_{\epsilon, \lambda}$
satisfies the Palais-Smale condition. To this end, we have to prove some
technical lemmas.

\begin{lemma}\label{doalves}
If $(v_{n})$ be a Palais-Smale sequence for $I_{\epsilon,\lambda}$
in $H_{\epsilon,\lambda}$ such that $v_{n}\rightharpoonup v$ in
$H_{\epsilon,\lambda}$ for some $v$ in $H_{\epsilon,\lambda}$, then
$$
I_{\epsilon,\lambda}(\widetilde{v}_{n})=I_{\epsilon,\lambda}(v_{n})-I_{\epsilon,\lambda}(v)+o_{n}(1)
$$
and
$$
I'_{\epsilon,\lambda}(\widetilde{v}_{n})=o_{n}(1)
$$
where $\widetilde{v}_{n}=v_{n}-v$.
\end{lemma}
\noindent {\bf Proof.} \, Firstly, we observe that the limits below
hold
\begin{equation} \label{GE1}
\displaystyle\int_{\mathbb{R}^{N}}F(|\widetilde{v}_{n}|^{2})
=\displaystyle\int_{\mathbb{R}^{N}}F(|v_{n}|^{2})-\displaystyle\int_{\mathbb{R}^{N}}F(|v|^{2})+o_{n}(1)
\end{equation}
and
\begin{equation} \label{GE2}
\displaystyle\int_{\mathbb{R}^{N}}\biggl|f(|\widetilde{v}_{n}|^{2})\widetilde{v}_{n}
-f(|v_{n}^{2}|)v_{n} +f(|v|^{2})v\biggl|^{r}=o_{n}(1)
\end{equation}
for some $r \in (2,2^{*})$.

We will show only the first limit because the same arguments can be
used in the proof of other one. We begin remarking that
$$
F(|\widetilde{v}_{n}+v|^{2})-F(|\widetilde{v}_{n}|^{2})=\displaystyle\int^{1}_{0}\displaystyle\frac{d}{dt}F(|\widetilde{v}_{n}+
t v|^{2})dt.
$$
Then
$$
F(|v_{n}|^{2})-F(|\widetilde{v}_{n}^{2}|)=\displaystyle\int^{1}_{0}f(|\widetilde{v}_{n}+
t v|^{2})2|\widetilde{v}_{n}+ t v||v|dt.
$$
From this, for each $\gamma >0$, there exists $C_{\gamma}>0$ such
that,
$$
F(|v_{n}|^{2})-F(|\widetilde{v}_{n}^{2}|)\leq
\displaystyle\int^{1}_{0}\bigl[2\gamma|\widetilde{v}_{n}+ t
v||v|+2C_{\gamma}|\widetilde{v}_{n}+ t v|^{2^{*}-1}|v|\bigl] dt.
$$
Hence
$$
|F(|{v}_{n}|^{2}))-F(|\widetilde{v}_{n}|^{2})|\leq C\bigl[
\gamma|\widetilde{v}_{n}|^{2} + C_{\gamma}|v|^{2}+
\gamma|\widetilde{v}_{n}|^{2^{*}}+ C_{\gamma}|v|^{2^{*}}\bigl]
$$
and so,
$$
|F(|{v}_{n}|^{2})-F(|\widetilde{v}_{n}|^{2})- F(|v|^{2})|\leq
C\bigl[ \gamma|\widetilde{v}_{n}|^{2} + \widehat{C}_{\gamma}|v|^{2}+
\gamma|\widetilde{v}_{n}|^{2^{*}}+
\widehat{C}_{\gamma}|v|^{2^{*}}\bigl],
$$
for some positive constant  $\widehat{C}_{\gamma}>0$. Now, repeating
the same arguments found in  \cite{Alves}, it follows that
$$
\int_{\mathbb{R}^{N}}|F(|{v}_{n}|^{2})-F(|\widetilde{v}_{n}|^{2})-
F(|v|^{2})| \to 0
$$
or equivalently,
$$
\int_{\mathbb{R}^{N}}F(|\widetilde{v}_{n}|^{2})=\int_{\mathbb{R}^{N}}F(|{v}_{n}|^{2})-\int_{\mathbb{R}^{N}}F(|v|^{2})
+ o_{n}(1).
$$
On the other hand,
\begin{equation} \label{GE5}
\|\widetilde{v}_{n}\|^{2}_{\epsilon,\lambda}=
\|v_{n}\|^{2}_{\epsilon,\lambda}-\|v\|^{2}_{\epsilon,\lambda}+o_{n}(1).
\end{equation}
Now, using (\ref{GE1}), (\ref{GE2}) and (\ref{GE5}), we deduce that
$$
I_{\epsilon,\lambda}(\widetilde{v}_{n})=I_{\epsilon,\lambda}(v_{n})-I_{\epsilon,\lambda}(v)+o_{n}(1)
$$
and
$$
I'_{\epsilon,\lambda}(\widetilde{v}_{n})=o_{n}(1),
$$
which completes the proof.  \hfill\rule{2mm}{2mm}

\begin{lemma}
\label{L1} Suppose that $f$ satisfies $(f_{1})-(f_{3})$. Let
$(v_{n}) \subset H_{\epsilon,\lambda}$ be a $(PS)_c$ sequence for
$I_{\epsilon,\lambda}$. Then there exists a constant $K>0$,
independent of $\epsilon$ and $\lambda$, such that
$$
\limsup_{n \to \infty}\|v_{n}\|^{2}_{\epsilon,\lambda} \leq K,
$$
for all $\epsilon, \lambda>0$.
\end{lemma}
\noindent {\bf Proof.} \, By  $(f_3)$,
\begin{eqnarray*}
c+o_{n}(1)\|v_{n}\|_{\epsilon,\lambda}&=&
I_{\epsilon,\lambda}(v_{n})-\frac{1}{\theta}I'_{\epsilon,\lambda}(v_{n})v_{n}\\
&\geq &
\bigl(\frac{1}{2}-\frac{1}{\theta}\bigl)\|v_{n}\|^{2}_{\epsilon,\lambda},
\end{eqnarray*}
where $o_{n}(1) \to 0$ as $n \to \infty$. Thus, we conclude that $(v_{n})$ is bounded in
$H_{\epsilon,\lambda}$ with
$$
\displaystyle\limsup_{n\rightarrow \infty}\|v_{n}\|^{2}_{\epsilon,
\lambda}\leq K:= \frac{2c\theta}{\theta-2}.
$$
\hfill\rule{2mm}{2mm}

\begin{lemma} \label{L2}
Suppose that $f$ satisfies $(f_{1})-(f_{3})$. Let $(v_{n}) \subset
H_{\epsilon,\lambda}$ be a $(PS)_c$ sequence for
$I_{\epsilon,\lambda}$. Then $c \geq 0$, and if $c=0$, we have that
$v_{n} \to 0$ in $H_{\epsilon,\lambda}$.
\end{lemma}
\noindent {\bf Proof}
As in the proof of Lemma \ref{L1},
\begin{equation} \label{E1}
c+o_{n}(1)\|v_{n}\|_{\epsilon,\lambda}=I_{\epsilon,\lambda}(v_{n})-\frac{1}{\theta}I'_{\epsilon,\lambda}(v_{n})v_{n}
\geq
\Big(\frac{1}{2}-\frac{1}{\theta}\Big)\|v_{n}\|^{2}_{\epsilon,\lambda}\geq
0
\end{equation}
that is
$$
c+o_{n}(1)\|v_{n}\|_{\epsilon,\lambda}\geq 0.
$$
The boundedness of $(v_{n})$ in $H_{\epsilon,\lambda}$ gives $c \geq
0$ after passage to the limit as $n \to \infty$. If $c=0$, the
inequality (\ref{E1}) gives  $v_{n} \to 0$ in $X_{
\epsilon,\lambda}$ as $n \to \infty$, finishing the proof of Lemma
\ref{L2}. \hfill\rule{2mm}{2mm}

\begin{lemma} \label{L3}
Suppose that $f$ satisfies $(f_{1})-(f_{3})$. Let $c>0$ and
$(v_{n})$ be a $(PS)_{c}$ sequence for $I_{\epsilon,\lambda}$. Then,
there exists $\delta>0$ such that
$$
\liminf_{n \to \infty}\int_{\mathbb{R}^{N}}|v_{n}|^{q} \geq \delta,
$$
with $\delta$ being independent of $\lambda$ and $\epsilon$.
\end{lemma}
\noindent {\bf Proof} From $(f_{1})-(f_{2})$, there exists a
constant $C>0$ such that
\begin{equation} \label{E2}
|f(t)t|\leq \frac{1}{4}|t|^{2}+C|t|^{q}
\end{equation}
for all $t \in \mathbb{R}$. Now, combining (\ref{E2}) with $I'_{\epsilon,\lambda}(v_{n})v_{n}=o_{n}(1)$, we get

\begin{equation} \label{E3}
\frac{1}{2}\|v_{n}\|^{2}_{\epsilon,\lambda}\leq
C\displaystyle\int_{\mathbb{{R}^{N}}}|v_{n}|^{q} + o_{n}(1).
\end{equation}

\noindent Recalling that
$$
\frac{1}{2}\|v_{n}\|^{2}_{\epsilon,\lambda}=I_{
\epsilon,\lambda}(v_{n})+\int_{\mathbb{R}^{N}}F(|v_{n}|^{2}),
$$
$(f_{3})$ combined with $I_{\epsilon,\lambda}(v_{n})=c+o_{n}(1)$ yields
\begin{equation} \label{E5}
\liminf_{n \to \infty} \|v_{n}\|^{2}_{\epsilon,\lambda} \geq 2c>0.
\end{equation}
Hence, the lemma follows from (\ref{E3}) and (\ref{E5}).
\hfill\rule{2mm}{2mm}

\begin{lemma} \label{L4}
Suppose that $f$ satisfies $(f_{1})-(f_{3})$ and $W$ satisfies
$(W_{1}')-(W_{2})$. Let $d>0$ be an arbitrary number. Given any
$\epsilon>0$ and $\eta>0$, there exist $\Lambda_{\eta}>0$ and
$R_{\eta}>0$, which are independent of $\epsilon$, such that if
$(v_{n})$ is a $(PS)_{c}$ sequence for $I_{\epsilon,\lambda}$ with
$c \leq d$ and $\lambda \geq \Lambda_{\eta}$, then
$$
\limsup_{n \to \infty}\int_{\mathbb{R}^{N} \setminus
B_{R_{\eta}}(0)}|v_{n}|^{q}<\eta.
$$
\end{lemma}
\noindent {\bf Proof.} \, Given any $R>0$, define
$$
X(R)=\{x \in \mathbb{R}^{N}: \, |x|>R; \, W(\epsilon x)\geq K_{0}\}
$$
and
$$
Y(R)=\{x \in \mathbb{R}^{N}:\, |x|>R; \, W(\epsilon x)<K_{0}\}.
$$

\noindent  Observe that
$$
\int_{X(R)}|v_{n}|^{2}\leq \frac{1}{\lambda
K_{0}+1}\int_{X(R)}(\lambda W(\epsilon x)+1)|v_{n}|^{2} \leq \|v_{n}\|^{2}_{\epsilon,\lambda}.
$$

\noindent From Lemma \ref{L1}, there exists $K>0$ such that
\begin{equation} \label{E6}
\limsup_{n \to \infty}\int_{X(R)}|v_{n}|^{2}\leq \frac{K}{\lambda
K_{0}+1}.
\end{equation}

\noindent On the other hand, by H\"{o}lder inequality
$$
\int_{Y(R)}|v_{n}|^{2} \leq
\Big(\int_{Y(R)}|v_{n}|^{2^{*}}\Big)^{\frac{2}{2^{*}}}(\mu(Y(R)))^{\frac{2}{N}}.
$$
Using Sobolev Embedding Theorem and Lemma \ref{L1}, there exists a
constant $\widehat{K}>0$ such that
\begin{equation} \label{E7}
\limsup_{n \to \infty}\int_{Y(R)}|v_{n}|^{2} \leq \widehat{K}
(\mu(Y(R)))^{\frac{2}{N}},
\end{equation}
where the constant $\widehat{K}$ is uniform on $c \in [0,d]$. Since
$$
Y(R) \subset \{x \in \mathbb{R}^{N}:\, W(\epsilon x) \leq K_{0}\},
$$
it follows from $(W_{2})$
\begin{equation}\label{E8}
\lim_{R \to \infty}\mu(Y(R))=0.
\end{equation}
Using an interpolation property, we derive that
$$
|v_{n}|_{L^{q}(\mathbb{R}^{N} \setminus B_{R}(0))} \leq
|v_{n}|^{\alpha}_{L^{2}(\mathbb{R}^{N} \setminus
B_{R}(0))}|v_{n}|^{1-\alpha}_{L^{2^{*}}(\mathbb{R}^{N} \setminus
B_{R}(0))}
$$
for some $\alpha \in (0,1)$, and thus by Lemma \ref{L1}, there
exists a constant $\widetilde{K}>0$ such that
\begin{equation}\label{E9}
\limsup_{n \to \infty}\int_{\mathbb{R}^{N} \setminus
B_{R}(0)}|v_{n}|^{q} \leq \widetilde{K }\limsup_{n \to
\infty}\Big(\int_{\mathbb{R}^{N} \setminus
B_{R}(0)}|v_{n}|^{2}\Big)^{\frac{q\alpha}{2}}.
\end{equation}
Combining (\ref{E6}) with (\ref{E7}) and (\ref{E8}), given $\eta>0$,
we can fix $R=R_{\eta}$ and $\Lambda_{\eta}>0$ such that
\begin{equation}\label{E10}
\limsup_{n \to \infty}\int_{\mathbb{R}^{N} \setminus
B_{R}(0)}|v_{n}|^{2} \leq
\Big(\frac{\eta}{2\widetilde{K}}\Big)^{\frac{2}{q\alpha}}
\end{equation}
for all $\lambda \geq \Lambda_{\eta}$. Consequently, from (\ref{E9})
and (\ref{E10}),
$$
\limsup_{n \to \infty} \int_{\mathbb{R}^{N} \setminus
B_{R}(0)}|v_{n}|^{q} \leq \eta.
$$
This concludes the proof of the lemma. \hfill\rule{2mm}{2mm}

\vspace{0.5 cm}

As a first consequence of the last lemma, we have the following result

\begin{corollary}\label{C1}
If $(v_{n})$ is a $(PS)_c$ sequence for $I_{\epsilon,\lambda}$ and
$\lambda$ is large enough, then its weak limit is nontrivial
provided that $c>0$.
\end{corollary}

\vspace{0.3 cm}

The next result we will prove functional $I_{\epsilon,\lambda}$ satisfies the Palais-Smale condition for
$\lambda$ sufficiently large for $\epsilon$ arbitrary. More
precisely, we state:

\vspace{0.3 cm}

\begin{proposition} \label{P1}
Suppose that $(f_{1})-(f_{3})$ and $(W_{1}')-(W_{2})$ hold. Then for
any $d>0$ and $\epsilon >0$, there exists $\Lambda>0$, independent of
$\epsilon$, such that $I_{ \epsilon,\lambda}$ satisfies the
$(PS)_{c}$ condition for all $c \leq d, \lambda \geq \Lambda$ and
$\epsilon
>0$. That is, any sequence $(v_{n}) \subset H_{\epsilon,\lambda}$
satisfying
\begin{equation}
I_{\epsilon,\lambda}(v_{n}) \to c \,\,\, \mbox{and} \,\,\, I'_{
\epsilon,\lambda}(v_{n}) \to 0, \label{ps}
\end{equation}
for $c \leq d$, has a strongly convergent subsequence in $H_{
\epsilon,\lambda}$.
\end{proposition}
\noindent {\bf Proof} Given any $d>0$ and $\epsilon>0$, take $c \leq
d$ and let $(v_{n})$ be a $(PS)_{c}$ sequence for
$I_{\epsilon,\lambda}$. From Lemma \ref{L1}, there are a subsequence
still denoted by $(v_{n})$ and $v \in H_{\epsilon,\lambda}$ such
that $(v_{n})$ is weakly convergent to $v$ in
$H_{\epsilon,\lambda}$. If $\widetilde{v}_{n}=v_{n}-v$, from Lemma
\ref{doalves},
\begin{equation} \label{E11}
I_{\epsilon,\lambda}(\widetilde{v}_{n})=I_{
\epsilon,\lambda}(v_{n})-I_{\epsilon,\lambda}(v)+o_{n}(1)
\end{equation}
and
\begin{equation} \label{E12}
I'_{\epsilon,\lambda}(\widetilde{v}_{n}) \to 0.
\end{equation}
Once that $I'_{ \epsilon,\lambda}(v)=0$, $(f_{3})$ gives

\begin{equation} \label{E13}
I_{ \epsilon,\lambda}(v)=I_{
\epsilon,\lambda}(v)-\frac{1}{\theta}I'_{ \epsilon,\lambda}(v)v \geq
\Big(\frac{1}{2}-\frac{1}{\theta}\Big)\|v\|^{2}_{\epsilon,\lambda,}\geq
0.
\end{equation}

Setting $c'=c-I_{\epsilon,\lambda}(v)$, by (\ref{E11})-(\ref{E13}),
we deduce that $c' \leq d$ and $(\widetilde{v}_{n})$ is a
$(PS)_{c'}$ sequence for $I_{ \epsilon,\lambda,}$, thus by Lemma
\ref{L2}, we have $c' \geq 0$. We claim that $c'=0$. On the
contrary, suppose that $c'>0$. From Lemma \ref{L3}, there is
$\delta>0$ such that
\begin{equation} \label{E14}
\liminf_{n \to
\infty}\int_{\mathbb{R}^{N}}|\widetilde{v}_{n}|^{q}>\delta.
\end{equation}
Letting $\eta=\frac{\delta}{2}$ and applying Lemma \ref{L4}, we get
$\Lambda >0$ and $R>0$ such that
\begin{equation} \label{E15}
\limsup_{n \to \infty}\int_{\mathbb{R}^{N} \setminus
B_{R}(0)}|\widetilde{v}_{n}|^{q} < \frac{\delta}{2}
\end{equation}
for the corresponding $(PS)_{c'}$ sequence for $I_{
\epsilon,\lambda}$ for all $\lambda \geq \Lambda$. Combining
(\ref{E14}) with (\ref{E15}) and using the fact that
$\widetilde{v}_{n} \rightharpoonup 0$ in $H_{ \epsilon,\lambda}$, we
derive
$$
\delta \leq \liminf_{n \to
\infty}\int_{\mathbb{R}^{N}}|\widetilde{v}_{n}|^{q} \leq \limsup_{n
\to \infty}\int_{\mathbb{R}^{N} \setminus
B_{R}(0)}|\widetilde{v}_{n}|^{q}\leq \frac{\delta}{2}
$$
which is impossible, then $c'=0$. Thereby,  by Lemma \ref{L2},
$\widetilde{v}_{n} \to 0$ in $H_{\epsilon,\lambda}$, that is, $v_{n}
\to v$ in $H_{\epsilon,\lambda}$ and the proof of Proposition
\ref{P1} is complete. \hfill\rule{2mm}{2mm}

\vspace{0,2cm}

In closing this section, we proceed with the study of
$(PS)_{c,\infty}$ sequences, that is, sequences $(v_{n})$ in
$H_{\epsilon,\lambda}$ verifying:
$$
\begin{array}{l}
i)\,\,\,\,\,\, \lambda_{n} \to \infty \\
\mbox{}\\
ii)\,\,\,\,\,\, (I_{\epsilon,\lambda_{n}}(v_{n})) \,\, \mbox{is bounded} \\
\mbox{}\\
iii)\,\,\,\,\,\, \|I'_{\epsilon,\lambda_{n}}(v_{n})\|^{*}_{\epsilon,\lambda_{n}} \to 0\\
\end{array}
$$
where $\|\,\,\,\|^{*}_{\epsilon,\lambda_{n}}$ is defined by
$$
\|\varphi\|^{*}_{\epsilon,\lambda_{n}}= \sup \{|\varphi(u)|;\  u \in
H_{\epsilon,\lambda_{n}},\,  \|u\|_{\epsilon,\lambda}\leq 1\}
 \,\,\,\, \mbox{for} \,\,\,\, \varphi \in H^{*}_{\epsilon,\lambda_{n}}.
$$

\begin{proposition}\label{P2}
Suppose that $(f_{1})-(f_{3})$ and $(W_{1}')-(W_{2})$ hold. Assume
that $(v_{n}) \subset H^1(\mathbb{R}^N,\mathbb{C})$ is a
$(PS)_{c,\infty}$ sequence. Then for each $\epsilon >0$ fixed, there
exists a subsequence still denoted by $(v_{n})$ and $v_{\epsilon}
\in H^1(\mathbb{R}^N,\mathbb{C})$ such that
\begin{itemize}
\item[i)]\ $v_{n} \to v_{\epsilon}$ in $H^1(\mathbb{R}^N,\mathbb{C})$. Moreover, $v_{\epsilon}=0$ on $\Omega_{\epsilon}^{c}$ and
$v_{\epsilon}$ is a solution of
$$
\left\{
 \begin{array}{ll}
\left( \displaystyle\frac{1}{i}\nabla - A(\epsilon z)\right)^2 u + u =
f(|u|^2)u,&~~z \in \Omega_{\epsilon}, \ \
\\
u \in H^1_{0}(\Omega_{\epsilon},\mathbb{C}),
\end{array}
       \right.
$$
where $\Omega_{\epsilon}=\frac{\Omega}{\epsilon}$.

\item[ii)]\ $\displaystyle
\lambda_{n}\int_{\mathbb{R}^{N}}W(\epsilon x)|v_{n}|^{2} \to 0$.

\item[iii)]\ $\displaystyle \|v_{n}-v\|^{2}_{\epsilon,\lambda_{n}}
\to 0$.
\end{itemize}
\end{proposition}

\vspace{0.2 cm}

\noindent \textbf{Proof.} As in the proof of Lemma \ref{L1}, the
sequence $(\|v_{n}\|_{\epsilon,\lambda_{n}})$ is bounded in
$\mathbb{R}$. Thus, we can extract a subsequence $v_{n}
\rightharpoonup v_{\epsilon}$ weakly in $H_{\epsilon,\lambda}$. For each
$m \in \mathbb{N}$, we fix
$$
C_{m}=\Big\{x \in \mathbb{R}^{N}:\, W_{\epsilon}(x)\geq
\frac{1}{m}\Big\},
$$
where $W_{\epsilon}(x)=W(\epsilon x).$ Hence, for each $n \in \mathbb{N}$,
$$
\int_{C_{m}}|v_{n}|^{2} \leq
m\int_{C_{m}}W_{\epsilon}(x)|v_{n}|^{2}\leq
\frac{m}{\lambda_{n}}\int_{\mathbb{R}^{N}}(1+\lambda_{n}W_{\epsilon}(x))|v_{n}|^{2}\leq
\frac{m}{\lambda_{n}}\|v_{n}\|^{2}_{\epsilon,\lambda_{n}}.
$$
Then, by Lemma \ref{L1},
$$
\int_{C_{m}}|v_{n}|^{2} \leq \frac{mK}{\lambda_{n}} \,\,\,
\mbox{for} \,\,\, n \in \mathbb{N},
$$
for some constant $K>0$. Using Fatou's Lemma, we get
$$
\int_{C_{m}}|v_{\epsilon}|^{2}=0
$$
after to passage to the limit as $n \to \infty$. Thus
$v_{\epsilon}=0$ almost everywhere in $C_{m}$. Observing that
$$
\mathbb{R}^{N} \setminus
W_{\epsilon}^{-1}(0)=\cup_{m=1}^{\infty}C_{m},
$$
we deduce that $v_{\epsilon}=0$ almost everywhere in $\mathbb{R}^{N}
\setminus W_{\epsilon}^{-1}(0)$.
 Now, recalling that $W_{\epsilon}^{-1}(0)=\overline{\Omega}_{\epsilon} \cup D_{\epsilon}$ and
 $\mu(D_{\epsilon})= \mu(\frac{1}{\epsilon}D)=0$, it follows that $v_{\epsilon}=0$ almost everywhere
 in
 $\mathbb{R}^{N} \setminus \overline{\Omega}_{\epsilon}$. As $\partial \Omega_{\epsilon}$ is a smooth set, let us conclude that
 $v_{\epsilon} \in H^{1}_{0}(\Omega_{\epsilon},\mathbb{C})$.

Arguing as in Lemma \ref{L4}, we can assert that given $\eta>0$
there exists $R>0$ such that
\begin{equation}
\limsup_{n \to \infty}\int_{\mathbb{R}^{N} \setminus
B_{R}(0)}|v_{n}|^{2} < \eta \label{lim1}
\end{equation}
 and
\begin{equation}
\limsup_{n \to \infty}\int_{\mathbb{R}^{N} \setminus
B_{R}(0)}|v_{n}|^{q} < \eta. \label{lim2}
\end{equation}
From $(f_{1})-(f_{2})$, for each $\gamma >0$ there exists
$C_{\gamma}>0$ such that
$$
|f(s)|\leq \gamma +C_{\tau}|s|^{(q-2)/2} \,\,\, \mbox{for all}
\,\,\, s \in \mathbb{R}.
$$
This inequality combined with Sobolev Embedding Theorems and the
limits ({\ref{lim1}) and (\ref{lim2}) yields  there is  a
subsequence, still denoted by $(v_{n})$, such that
\begin{equation} \label{E16}
\lim_{n \to
\infty}\int_{\mathbb{R}^{N}}f(|v_{n}|^{2})v_{n}=\int_{\mathbb{R}^{N}}f(|v_{\epsilon}|^{2})v_{\epsilon}
\end{equation}
and
\begin{equation} \label{E17}
\lim_{n \to
\infty}\int_{\mathbb{R}^{N}}f(|v_{n}|^{2})v_{\epsilon}=\int_{\mathbb{R}^{N}}f(|v_{\epsilon}|^{2})v_{\epsilon}.
\end{equation}
Thus, by (\ref{E16}) and (\ref{E17}),
\begin{eqnarray*}
\|v_{n}-v_{\epsilon}\|^{2}_{\epsilon}\leq
\|v_{n}-v_{\epsilon}\|^{2}_{\epsilon,\lambda_{n}}&=&
\|v_{n}\|^{2}_{\epsilon,\lambda_{n}}-
\|v_{\epsilon}\|^{2}_{\epsilon,\lambda_{n}}+o_{n}(1)\\
&=& I'_{\epsilon,\lambda_{n}}(v_{n})v_{n}-
I'_{\epsilon,\lambda_{n}}(v_{\epsilon})v_{\epsilon}+o_{n}(1)=o_{n}(1).
\end{eqnarray*}
Thereby, $v_{\epsilon}$ is a solution of
$$
\left\{
 \begin{array}{ll}
\left( \displaystyle\frac{1}{i}\nabla - A(\epsilon z)\right)^2 u + u =
f(|u|^2)u,&~~z \in \Omega_{\epsilon}, \ \
\\
u \in H^1_{0}(\Omega_{\epsilon},\mathbb{C}),
\end{array}
       \right. \eqno{(D_{\epsilon})}
$$
and the proof of $i)$ is complete.

To deduce $ii)$, we start observing that
$$
\|v_{n}\|^{2}_{\epsilon} +
\lambda_{n}\displaystyle\int_{\mathbb{R}^{N}}W_{\epsilon}(x)|v_{n}|^{2}=\displaystyle\int_{\mathbb{R}^{N}}f(|v_{n}|^{2})|v_{n}|^{2}+o_{n}(1).
$$
Since $v_{n}\rightarrow v_{\epsilon}$ in
$H^{1}(\mathbb{R}^{N},\mathbb{C})$ and $v_{\epsilon}$ is a solution
of $(D_{\epsilon})$, we obtain
$$
\displaystyle\lim_{n \to \infty}
\lambda_{n}\int_{\mathbb{R}^{N}}W_{\epsilon}(x)|v_{n}|^{2}=0,
$$
which we conclude ii).

For to prove iii), we observe that
$$
\|v_{n}-v_{\epsilon}\|^{2}_{\epsilon,\lambda_{n}}=
\|v_{n}-v_{\epsilon}\|^{2}_{\epsilon} +
\lambda_{n}\displaystyle\int_{\mathbb{R}^{N}}W_{\epsilon}(x)|v_{n}|^{2}.
$$
This last equality allow that we conclude that
$$
\|v_{n}-v_{\epsilon}\|^{2}_{\epsilon,\lambda_{n}}=o_{n}(1).
$$
 $\hfill \rule{2mm}{2mm}$

\begin{corollary} \label{T2}
Suppose that $(W_{1}')-(W_{2})$ and $(f_{1})-(f_{4})$ hold. Then for
each $\epsilon >0$ and a sequence $(v_{n})$ of solutions of
$(P_{\epsilon,\lambda_{n}})$ with $\lambda_{n} \to \infty$ and
$\displaystyle\limsup_{n\to \infty}I_{\epsilon,\lambda_n}(v_{n})< \infty$, there
exists a subsequence that converges strongly in
$H^{1}(\mathbb{R}^{N},\mathbb{C})$ to a solution of the problem
$(D_{\epsilon})$.
\end{corollary}

\noindent \textbf{Proof.} By assumptions, there exist $c\in
\mathbb{R}$ and a subsequence of $(v_n)$, still denoted by $(v_n)$,
such that $(v_n)$ is a $(PS)_{c,\infty}$ sequence. The rest of the
proof follows from Proposition \ref{P2}. $\hfill \rule{2mm}{2mm}$


\section{Behavior of minimax levels}

This section is devoted to the study of the behavior of the minimax
levels with respect to parameter $\lambda$ and $\epsilon$. For this
purpose, we introduce some notations. In the next,
$\mathcal{M}_{\epsilon,\lambda}$ denotes the Nehari manifold
associated to $I_{\epsilon,\lambda}$, that is,
$$
\mathcal{M}_{\epsilon,\lambda}=\Big\{v \in H_{\epsilon,\lambda}:\,
v\not=0 \,\,\, \mbox{and} \,\,\, I'_{\epsilon,\lambda}(v)v=0\Big\}
$$
and
$$
c_{\epsilon,\lambda}=\inf_{v \in
\mathcal{M}_{\epsilon,\lambda}}I_{\epsilon,\lambda}(v).
$$
From $(f_1)-(f_4)$ and arguing as in \cite{Willem}, we can prove
that $c_{\epsilon,\lambda}$ can also characterized as the mountain
pass minimax value associated with $I_{\epsilon,\lambda}$.

On account of the proof of Proposition \ref{P2}, when $\lambda$ is
large, problem $(D_{\epsilon})$ can be seen as a limit problem of
$(D_{\epsilon,\lambda,})$ for each $\epsilon>0$. The functional
corresponding to the problem $(D_{\epsilon})$ is given by
$$
P_{\epsilon}(v)=\frac{1}{2}\int_{\Omega_{\epsilon}}(|\nabla_{\epsilon}
v|^{2}+|v|^{2})-\frac{1}{2}\int_{\Omega_{\epsilon}}F(|v|^{2})
$$
for every $v \in H_{0}^{1}(\Omega_{\epsilon},\mathbb{C})$. Here and
subsequently, $\mathcal{M}_{\epsilon}$ denotes the Nehari manifold
associated to $P_{\epsilon}$ and
$$
c(\epsilon,\Omega)=\displaystyle \inf_{v \in
\mathcal{M}_{\epsilon}}P_{\epsilon}(v)
$$
stands for the mountain pass minimax associated with $P_{\epsilon}$.

Since $0 \in \Omega$, there is $r>0$ such that $B_{r}=B_{r}(0)
\subset \Omega$ and
$B_{\frac{r}{\epsilon}}=B_{\frac{r}{\epsilon}}(0) \subset
\Omega_{\epsilon}$. We will denote by
$P_{\epsilon,B_{r}}:H_{0}^{1}(B_{\frac{r}{\epsilon}}(0),\mathbb{C})
\to \mathbb{R}$ the functional
$$
P_{\epsilon,B_{r}}(v)=\frac{1}{p}\int_{B_{\frac{r}{\epsilon}}}(|\nabla_{\epsilon}
v|^{2}+|v|^{2})-\int_{B_{\frac{r}{\epsilon}}}F(|v|^{2}).
$$
Furthermore, we write $\mathcal{M}_{\epsilon,B_{r}}$ the Nehari
manifold associated to $P_{\epsilon,B_{r}}$ and
$$
c(\epsilon,B_{r})=\displaystyle \inf_{v \in
\mathcal{M}_{\epsilon,B_{r}}}P_{\epsilon,B_r}(v).
$$

The next Lemma will be useful for studying the behavior of
$c_{\epsilon,\lambda}$ as $\lambda$ goes to infinity.

Hereafter, we denote by $b_{\infty}$ the mountain pass level associated with functional $J: H^{1}(\mathbb{R}^{N},\mathbb{R}) \to \mathbb{R}$, given by
$$
J_{\infty}(v)=\frac{1}{2}\displaystyle\int_{\mathbb{R}^{N}}|\nabla
v|^{2}+\frac{1}{2}\displaystyle\int_{\mathbb{R}^{N}}|v|^{2}-\frac{1}{2}\displaystyle\int_{\mathbb{R}^{N}}F(|v|^{2}).
$$
From diamagnetic inequality (\ref{diamagnetica}),
$$
J_{\infty}(tv)\leq I_{\epsilon,\lambda}(tv) \,\, \forall t \geq 0,
$$
and so,
$$
b_{\infty} \leq c_{\epsilon,\lambda}.
$$

The following result is a consequence of Proposition \ref{P2}.

\begin{proposition} \label{P3}
Suppose $(f_{1})-(f_{4})$ and $(W_{1}')-(W_{2})$ hold. Let $\epsilon
>0$ be an arbitrary number. Then,
$$
\lim_{\lambda \to \infty}c_{\epsilon,\lambda}=c(\epsilon,\Omega).
$$
\end{proposition}

\noindent \textbf{Proof.} By Proposition \ref{P1} and Mountain Pass
Theorem, we can assume that there are two sequences,  $\lambda_{n}
\to \infty$ and  $(v_{n}) \subset H_{\epsilon,\lambda_{n}}$,   such
that
$$
I_{\epsilon,\lambda_{n}}(v_{n})=c_{\epsilon,\lambda_{n}}>0\,\,\,
\mbox{and} \,\,\, I'_{\epsilon,\lambda_{n}}(v_{n})=0.
$$
From definitions of $c_{\epsilon,\lambda_{n}}$ and $c(\epsilon,
\Omega)$,
$$
c_{\epsilon,\lambda_{n}} \leq c(\epsilon, \Omega) \,\,\, \mbox{for
all} \,\,\, n \in \mathbb{N}
$$
which implies
$$
0\leq I_{\epsilon,\lambda_{n}}(v_{n}) \leq c(\epsilon, \Omega)
\,\,\, \mbox{and} \,\,\, I'_{\epsilon,\lambda_{n}}(v_{n})=0.
$$
Thus, for some subsequence $(v_{n_{j}})$, there exists $c \in
[0,c(\epsilon, \Omega)]$ such that
$$
I_{\epsilon,\lambda_{n_{j}}}(v_{n_{j}})=c_{\epsilon,\lambda_{n_{j}}}
\to c \,\,\, \mbox{and} \,\,\,
I_{\epsilon,\lambda_{n_{j}}}'(v_{n_{j}}) \to 0
$$
showing that $(v_{n_{j}})$ is a $(PS)_{c,\infty}$, and so,
$$
\displaystyle\int_{\mathbb{R}^{N}}|\nabla_{\epsilon} v_{n}|^{2}
+\displaystyle\int_{\mathbb{R}^{N}}(\lambda_{n}W_{\epsilon}(x)+1)|v_{n}|^{2}
\geq 2 c_{\epsilon, \lambda_ n} \geq 2 b_{\infty} >0 \,\,\, \forall n
\in \mathbb{N}.
$$
By Proposition \ref{P2},
$$
\lambda_{n}\int_{\mathbb{R}^{N}}W_{\epsilon}(x)|v_{n}|^{2} \to 0
\,\,\, \mbox{as} \,\,\, n \to +\infty
$$
then,
\begin{equation}\label{bom}
\displaystyle\int_{\mathbb{R}^{N}}|\nabla_{\epsilon} v_{n}|^{2}
+\displaystyle\int_{\mathbb{R}^{N}}|v_{n}|^{2}
 \geq 2 b_{\infty} >0 +o_n(1) \,\,\, \forall n \in \mathbb{N},
\end{equation}
implying that any subsequence of $(v_{n})$ does not converge to zero
in $H^{1}(\mathbb{R}^{N},\mathbb{C})$.

From Proposition \ref{P2}, there exist a subsequence
$(v_{n_{j_{k}}})$ and $v \in H^{1}(\mathbb{R}^{N},\mathbb{C})$ such
that
\begin{equation} \label{E24}
v_{n_{j_{k}}} \to v \,\,\, \mbox{strongly in} \,\,\,
H^{1}(\mathbb{R}^{N},\mathbb{C}) \,\,\, \mbox{and} \,\,\, v=0 \,\,\,
\mbox{in} \,\,\, \mathbb{R}^{N} \setminus \Omega_{\epsilon}.
\end{equation}
From (\ref{bom}) and (\ref{E24}), $v \not=0$ in
$H^{1}_{0}(\Omega_{\epsilon},\mathbb{C})$ and $v$ is a solution of the
problem $(D_{\epsilon})$, from where it follows that
\begin{equation} \label{E25}
P_{\epsilon}(v) \geq c(\epsilon,\Omega).
\end{equation}
On the other hand,
\begin{equation}\label{E26}
P_{\epsilon}(v)=\lim_{k\to
\infty}I_{\epsilon,\lambda_{n_{j_{k}}}}(v_{n_{j_{k}}})=\lim_{k\to
\infty}c_{\epsilon,\lambda_{n_{j_{k}}}}=c\leq c(\epsilon,\Omega).
\end{equation}
Therefore, (\ref{E25}) and (\ref{E26}) give
$$
\lim_{k\to
\infty}c_{\epsilon,\lambda_{n_{j_{k}}}}=c(\epsilon,\Omega).
$$
As a result,  $c_{\epsilon,\lambda} \to c(\epsilon,\Omega)$
as $\lambda \to \infty$, and the lemma follows. $\hfill\rule{2mm}{2mm}$

\vspace{0,2cm}

\begin{corollary} \label{T3}
Suppose that $(W_{1}')-(W_{2})$ and $(f_{1})-(f_{4})$ hold. Then for
each $\epsilon >0$ and a sequence $(v_{n})$ of least energy
solutions of $(D_{\epsilon,\lambda_{n}})$ with $\lambda_{n} \to
\infty$ and $\displaystyle\limsup_{n\to
\infty}I_{\epsilon,\lambda_n}(v_{n})< \infty$, there exists a
subsequence that converges strongly in
$H^{1}(\mathbb{R}^{N},\mathbb{C})$ to a least energy solution of the
problem $(D_{\epsilon})$.
\end{corollary}
\textbf{Proof.} The proof is a consequence of Propositions \ref{P2}
and \ref{P3}. $\hfill \rule{2mm}{2mm}$

\vspace{0.3 cm}

Hereafter, $r>0$ denotes a number such that $B_{r}(0)
\subset \Omega$ and the sets
$$
\Omega_{+}=\{x \in \mathbb{R}^{N}:\, d(x,\overline{\Omega})\leq r\}
$$
and
$$
\Omega_{-}=\{x \in \mathbb{R}^{N}:\, d(x,\partial \Omega) \geq r\}
$$
are homotopically equivalent to $\Omega$. The existence of this $r$
is given by condition $(W_{1}')$. For each $v \in
H^{1}(\mathbb{R}^{N},\mathbb{C})$ with compact support, we consider
the barycenter of $v$
$$
\beta(v)= \frac{\displaystyle \int_{\mathbb{R}^{N}}x
|v|^{2}}{\displaystyle \int_{\mathbb{R}^{N}}|v|^{2}}.
$$

Consider $R>0$ such that $\Omega \subset B_{R}(0)$  and set
the auxiliary function
$$
\xi_{\epsilon}(t)= \left\{
\begin{array}{l}
1, \,\, 0\leq t \leq \frac{R}{\epsilon}\\
\mbox{}\\
\frac{R}{\epsilon t},\,\,\,\, \frac{R}{\epsilon}\leq t.
\end{array}
\right.
$$
For $v \in H^{1}(\mathbb{R}^{N},\mathbb{C}),v\not=0$, define
$$
\beta_{\epsilon}(v)=\frac{\displaystyle
\int_{\mathbb{R}^{N}}x\xi_{\epsilon}(|x|)|v|^{2}}{\displaystyle
\int_{\mathbb{R}^{N}}|v|^{2}}.
$$

Now for each $y \in \mathbb{R}^{N}$ and $R>2{\rm{diam}}(\Omega)$
fix
$$
Z_{\frac{R}{\epsilon},\frac{r}{\epsilon},y}=\Big\{x \in
\mathbb{R}^{N}:\, \frac{r}{\epsilon}\leq |x-y|\leq
\frac{R}{\epsilon}\Big\}.
$$
We observe that if $y \notin \frac{1}{\epsilon}\Omega_{+}$, then
$\overline{\Omega_{\epsilon}} \cap
B_{\frac{r}{\epsilon}}(y)=\emptyset$. As a consequence
\begin{equation} \label{E27}
\overline{\Omega_{\epsilon}} \subset
Z_{\frac{R}{\epsilon},\frac{r}{\epsilon},y}
\end{equation}
for every $y \notin \frac{1}{\epsilon}\Omega_{+}$. Moreover, for $y
\in \mathbb{R}^{N}$, we will consider the number  $\alpha(R,r,\epsilon,y)$ given by
$$
\alpha(R,r,\epsilon,y)=\inf \Big\{\widehat{J}_{\epsilon,y}(v):\,
\beta(v)= y \,\,\, \mbox{and} \,\,\, v \in
\widehat{\mathcal{N}}_{\epsilon,y}\Big\}
$$
where
$$
\widehat{J}_{\epsilon,y}(v)=\frac{1}{2}\int_{Z_{\frac{R}{\epsilon},\frac{r}{\epsilon},y}}(|\nabla
v|^{2}+|v|^{2})-\frac{1}{2}\int_{Z_{\frac{R}{\epsilon},\frac{r}{\epsilon},y}}F(|v|^{2})
$$
and
$$
\widehat{\mathcal{N}}_{\epsilon,y}=\Big\{v \in
H^{1}_{0}(Z_{\frac{R}{\epsilon},\frac{r}{\epsilon},y},\mathbb{R}):v\not=0
\,\,\, \mbox{and} \,\,\, \widehat{J}_{\epsilon,y}'(v)v=0\Big\}.
$$
From now on, we will write $\alpha(R,r,\epsilon,0)$ as
$\alpha(R,r,\epsilon)$, $\widehat{J}_{\epsilon,0}$ as
$\widehat{J}_{\epsilon}$ and  $\widehat{\mathcal{N}}_{\epsilon,0}$
as $\widehat{\mathcal{N}}_{\epsilon}$.

\begin{lemma} \label{L6}
Assume that $(f_{1})-(f_{4})$ hold. Then, there exist $\epsilon^{*}, \delta>0$
such that
$$
b_{\infty} + \delta < \alpha(R,r,\epsilon)
$$
for every $\epsilon \in (0, \epsilon^{*})$.
\end{lemma}
\noindent \textbf{Proof.} See proof in \cite[Proposition 4.1]{Alves10}.  $\hfill \rule{2mm}{2mm}$

\vspace{0.5 cm}

\section{Proof of Theorem \ref{T1}}

For $r>0$ and $\epsilon>0$, let $v_{r\epsilon} \in
H^{1}_{0}(B_{\frac{r}{\epsilon}}(0))$ be a nonnegative radially
symmetric function such that
$$
J_{\epsilon,B_{r}}(v_{r\epsilon})=b(\epsilon,B_{r})\,\,\, \mbox{and}
\,\,\, J'_{\epsilon,B_{r}}(v_{r\epsilon})=0,
$$
where
$$
J_{\epsilon,B_{r}}(v)=
\frac{1}{2}\displaystyle\int_{B_{\frac{r}{\epsilon}}(0)}|\nabla
u|^{2} + \frac{1}{2}\displaystyle\int_{B_{\frac{r}{\epsilon}}(0)}|
u|^{2}- \frac{1}{2}\displaystyle\int_{B_{\frac{r}{\epsilon}}(0)}F(|
u|^{2})
$$
whose existence is proved in \cite{Alves10}. For $r>0$ and $\epsilon
>0$, define \linebreak $\Psi_{\epsilon}: \frac{1}{\epsilon}\Omega_{-} \to
\mathcal{M}_{\epsilon,\lambda}$ by
$$
\Psi_{\epsilon}(y)(x)= \left\{
\begin{array}{l}
t_{\epsilon,y}e^{i \tau_{\epsilon,y}(x)}v_{r\epsilon}(|x-y|), \,\, x \in B_{\frac{r}{\epsilon}}(y)\\
\mbox{}\\
0, \,\,x \notin B_{\frac{r}{\epsilon}}(y),
\end{array}
\right.
$$
where $\tau_{\epsilon, y}(x) := \sum_{j=1}^N A_j(\epsilon y)x^j$ and $t_{\epsilon,y} \in (0,+\infty)$ is such that
$$
t_{\epsilon,y}\textrm{e}^{i\tau_{\epsilon,y}(.)}v_{r \epsilon}(|\cdot -y|) \in
\mathcal{M}_{\epsilon,\lambda}.
$$
It is immediate that $\beta_{\epsilon}(\Psi_{\epsilon}(y))=y$ for all $y \in
\frac{1}{\epsilon}\Omega_{-}$.

\begin{lemma} \label{NovoLema}
Uniformly for $y \in \frac{1}{\epsilon}\Omega_-$, there holds
$$
\lim_{\epsilon \to
0}P_{\epsilon}(\Psi_{\epsilon}(y))=b_{\infty}.
$$
\end{lemma}

\noindent\textbf{Proof.}  Given three sequences
$\epsilon_n \to 0, \lambda_n \to +\infty$ and $(y_n) \subset \frac{1}{\epsilon_n}\Omega_{-}$,
we will prove that
$$
P_{\epsilon_n}(\Psi_{\epsilon_n}(y_n)) \to b_{\infty} ~~~~ \mbox{as} ~~~~ n \to +\infty.
$$
Let $t_n := t_{\epsilon_n,y_n}$ and $v_n=v_{r \epsilon_{n}}$ be as in the definition of
$\Psi_{\epsilon_n}$. Using the diamagnetic inequality, we have
\begin{equation}\label{ZZ1}
b(\epsilon_n,B_{r}) \leq P_{\epsilon_n}(\Psi_{\epsilon_n}(y_n)).
\end{equation}
On the other hand,
$$
P_{\epsilon_n}(\Psi_{\epsilon_n}(y_n)) \leq b(\epsilon_n,B_{r}) + \dfrac{t_{n}^2}{2}
\dis\int_{B_{\frac{r}{\epsilon_n}}(y_n)} |A(\epsilon_n y_n)-A(\epsilon_n x + \epsilon_x y_n)||v_n|^2\;dx
$$
A direct computation implies that $(t_n)$ is bounded, hence
\begin{equation}\label{ZZ2}
P_{\epsilon_n}(\Psi_{\epsilon_n}(y_n)) \leq b(\epsilon_n,B_{r}) + C_1 \dis\int_{\mathbb{R}^{N}} |A(\epsilon_n y_n)-A(\epsilon_n x + \epsilon_x y_n)||v_n|^2\;dx.
\end{equation}
Moreover, it is possible to prove that there is $v \in H^{1}(\mathbb{R}^{N})$ such that
$$
v_n \to v \,\,\ \mbox{in} \,\,\, H^{1}(\mathbb{R}^{N}).
$$
Once that $A$ is continuous and belongs to $L^{\infty}(\mathbb{R}^{N},\mathbb{R}^{N})$, the above limit limit loads to
\begin{equation} \label{ZZ3}
 \dis\int_{\mathbb{R}^{N}} |A(\epsilon_n y_n)-A(\epsilon_n x + \epsilon_x y_n)||v_n|^2\;dx.
\end{equation}
Combining (\ref{ZZ1}), (\ref{ZZ2})  and  (\ref{ZZ3}) with the limit $b(\epsilon_n,B_{r}) \to b_\infty$, we derive that
$$
P_{\epsilon_n}(\Psi_{\epsilon_n}(y_n)) \to b_\infty,
$$
finishing the proof. {\fim}

\vspace{0.5 cm}

In what follows, we set
\begin{equation} \label{def-g}
g(\epsilon) := \sup_{y \in \frac{1}{\epsilon}\Omega_{-}}|P_{\epsilon}(\Psi_{\epsilon}(y)) - b_{\infty}|.
\end{equation}
From Lemma \ref{NovoLema},  $g(\epsilon) \to 0$ as $\epsilon \to 0$. Consequently, we can assume that
\begin{equation} \label{Eg0}
g(\epsilon) < \delta \,\,\, \forall \epsilon \in (0, \epsilon^{*}),
\end{equation}
where $\epsilon^{*}, \delta$ were given in Lemma \ref{L6}. Here, we are assuming that we can decrease $\epsilon^{*}$ if necessary.

Moreover, we observe that $I_{\epsilon, \lambda}(\Psi_{\epsilon}(y))= P_{\epsilon}(\Psi_{\epsilon}(y)) \leq b_{\infty} + g(\epsilon)$ for all $y \in \frac{1}{\epsilon}\Omega_{-}$ and $\lambda >0$.
Hence, the set
$$
\mathcal{O}_{\epsilon, \lambda}:=\{ u \in \mathcal{M}_{\epsilon, \lambda}:
I_{\epsilon, \lambda}(u)\leq b_{\infty} +g(\epsilon)\}
$$
contains the functions $\Psi_{\epsilon}(y)$ for $y \in \frac{1}{\epsilon}\Omega_{-}$, showing that $\mathcal{O}_{\epsilon, \lambda} \not = \emptyset$.

\begin{lemma} \label{L7}
Suppose $(f_{1})-(f_{4})$ and $(W_{1}')-(W_{2})$ hold. Let
$\epsilon^{*}>0$ given by Lemma \ref{L6}. Then for any $\epsilon \in
(0, \epsilon^{*})$, there exists $\lambda^{*}>0$ which depends on
$\epsilon$ such that
$$
\beta_{\epsilon}(v) \in \frac{1}{\epsilon}\Omega_{+}
$$
for all $\lambda > \lambda^{*},0< \epsilon < \epsilon^{*}$ and $v
\in \mathcal{O}_{\epsilon,\lambda}$.
\end{lemma}
\noindent \textbf{Proof.} Suppose by contradiction that there exists
a sequence $(\lambda_{n})$ with $\lambda_{n} \to \infty$ such that
\begin{equation} \label{Eg1}
v_{n} \in \mathcal{M}_{\epsilon,\lambda_{n}}, \,\,
I_{\epsilon,\lambda_{n}}(v_{n}) \leq b_{\infty} + g(\epsilon)
\end{equation}
and
\begin{equation}\label{E30}
\beta_{\epsilon}(v_{n}) \notin \frac{1}{\epsilon}\Omega_{+}.
\end{equation}
Repeating the same arguments used in the proofs of Lemma \ref{L4}
and Proposition \ref{P2}, $(\|v_n\|_{\epsilon,\lambda_n})$ is a
bounded sequence in $\mathbb{R}$ and there exists $v \in
H^{1}(\mathbb{R}^N,\mathbb{C})$ such that $v_n \rightharpoonup v$
weakly in $H^{1}(\mathbb{R}^N,\mathbb{C})$, $v = 0$ in
$\mathbb{R}^N\setminus \Omega_\epsilon$ and for each $\eta>0$ there
exists $R>0$ such that
$$
\limsup_{n \to \infty}\int_{\mathbb{R}^{N} \setminus
B_{R}(0)}|v_{n}|^{2}<\eta.
$$
This fact implies that
\[
v_n \to v \,\,\, \mbox{strongly in} \,\,\,
L^p(\mathbb{R}^N,\mathbb{C}).
\]
Hence by interpolation,
\[
v_n \to v \,\,\, \mbox{strongly in} \,\,\,
L^t(\mathbb{R}^N,\mathbb{C})\,\,\, \mbox{for all} \,\,\, t\in
[2,2^*).
\]
On the other hand, since $v_n \in
\mathcal{M}_{\epsilon,\lambda_{n}}$,  from  (\ref{bom}),
$$
0<2b_{\infty} \leq
\displaystyle\int_{\mathbb{R}^{N}}f(|v_n|^{2})|v_n|^{2} + o_n(1), \,\,\, \mbox{for
all} \,\, n \in \mathbb{N},
$$
from where it follows that
$$
0< 2 b_{\infty} \leq \int_{\mathbb{R}^N}f(|v|^{2})|v|^{2},
$$
which yields
\begin{equation} \label{E31}
v \not=0, \, P'_{\epsilon}(v)v \leq 0 \,\,\,\, \mbox{and} \,\,\,\,
\lim_{n \to \infty}\beta_{\epsilon}(v_{n})=\beta(v).
\end{equation}
From (\ref{E30}) and (\ref{E31}), $y=\beta(v) \notin
\frac{1}{\epsilon}\Omega_{+}$, $\Omega_{\epsilon} \subset
A_{\frac{R}{\epsilon},\frac{r}{\epsilon},y}$ and there exists $\tau
\in (0,1]$ such that $\tau |v| \in
\widehat{\mathcal{N}}_{\epsilon,y}$. Thereby, combining diamagnetic
inequality (\ref{diamagnetica}) with (\ref{Eg0}) and (\ref{Eg1}), we get
$$
\widehat{J}_{\epsilon,y}(\tau |v|)\leq P_{\epsilon}(\tau v) \leq
\liminf_{n \to \infty}I_{\epsilon,\lambda_{n},}(\tau v_{n}) \leq
\liminf_{n \to \infty}I_{\epsilon,\lambda_{n}}(v_{n}) \leq
b_{\infty} + \delta
$$
which implies
$$
\alpha(R,r,\epsilon,y)\leq b_{\infty} + \delta.
$$
On the other hand, since
$$
\alpha(R,r,\epsilon,y)=\alpha(R,r,\epsilon)
$$
we have
$$
\alpha(R,r,\epsilon) \leq b_{\infty} + \delta,
$$
obtaining a contradiction with Lemma \ref{L6}, and the proof is complete. $\hfill
\rule{2mm}{2mm}$

\vspace{0.5 cm}

We claim that
\begin{equation}\label{E32}
cat \mathcal{O}_{\epsilon,\lambda} \geq cat (\Omega)
\end{equation}
for all $\epsilon \in (0, \epsilon^{*})$ and $\lambda \geq
\lambda^{*}$. In fact, suppose that
$$
\mathcal{O}_{\epsilon,\lambda}=\cup_{i=1}^{n}O_{i}
$$
where $O_{i},i=1,...,n$, is closed and contractible in
$\mathcal{O}_{\epsilon,\lambda}$, that is, there exists
$h_{i} \in C([0,1] \times
O_{i},\mathcal{O}_{\epsilon,\lambda})$ such that, for
every, $v \in O_{i}$,
$$
h_{i}(0,v)=v \,\,\, \mbox{and} \,\,\, h_{i}(1,u)=w_{i}
$$
for some $w_{i}\in \mathcal{O}_{\epsilon,\lambda}$.
Consider
$$
B_{i}=\Psi_{\epsilon}^{-1}(O_{i}), \,\,\, i=1,...,n.
$$
The sets $B_{i}$ are closed and
$$
\frac{1}{\epsilon}\Omega_{-}=B_{1}\cup...\cup B_{n}.
$$
Consider the deformation $g_{i}:[0,1] \times B_{i} \to
\frac{1}{\epsilon}\Omega_{+}$ given
$$
g_{i}(t,y)=\beta_{\epsilon}(h_{i}(t,\Psi_{\epsilon}(y))).
$$
From Lemma \ref{L7}, the function $g_{i}$ is well defined. Thus,
$B_{i}$ is contractile in $\frac{1}{\epsilon}\Omega_{+}$. Hence,
$$
cat(\Omega)= cat(\Omega_{\epsilon})=
 cat_{\frac{1}{\epsilon} \Omega_{+}}\big(\frac{1}{\epsilon}\Omega_{-}\big)\leq cat \, \mathcal{O}_{\epsilon,\lambda}
$$
which verifies (\ref{E32}).

Now, we are ready to conclude the proof of Theorem \ref{T1}. From
Proposition \ref{P1} the functional $I_{\epsilon,\lambda}$ satisfies
the Palais-Smale condition provided that $\lambda \geq \lambda^{*}$.
Thus, by Lusternik-Schirelman theory, the functional
$I_{\epsilon,\lambda}$ has at least $cat(\Omega)$ critical points
for all $\epsilon \in (0, \epsilon^{*})$ where $\epsilon^{*}>0$ is
given by Lemma \ref{L6}. The proof is complete. $\hfill
\rule{2mm}{2mm}$

\end{document}